\documentclass[12pt,a4paper]{article}
\usepackage{amsmath,amssymb,amsthm,mathptmx}
\newtheorem{theorem}{Theorem}

\author{Anastasia V. Parusnikova}
\title{On Gevrey orders \\of power expansions of solutions \\to the third and fifth 
Painlev\'{e} equations\thanks{This study was carried out within The National Research University Higher School of Economics Academic Fund Program in 2013-2014, research grant No. 12-01-0030 and within RFBR, Grant Nr. 12-01-31414-mol-a.}}

\date{}

\begin{document}
\maketitle

\begin{abstract}
The question under consideration is Gevrey summability of power expansions of solutions to the third and fifth Painlev\'{e} equations near infinity. Methods of French and Japaneese schools are used to analyse these properties of formal power-series solutions. The results obtained are compared with the ones obtained by means of Power Geometry.

\end{abstract}



\section{General theory}
~~~~Let $V$ be an open sector  with a vertex in the infinity on extended complex plane or on Riemann surface of logarithm, i.~e. $V=\{z:|z|>R,\mathrm{Arg}\,z\in(\varphi_1,\varphi_2)\}$. $w$ -is a holomorphic on $V$ function and $\hat{w}={\displaystyle{\sum_{k=0}^{\infty}a_k z^{-k}}}$ is some formal series belonging to $\mathbb{C}[[1/z]]$.

A function $w$ is said to be \textit{asymptotically approximated by a series $\hat{w}$ on $V$}, if for the points $z$ of any closed subsector $Y$ of $V$ and for any $n \in \mathbb{N}$ there exist the constants $M_{Y,n}>0$:
$$|z^{n}||w(z)-\sum_{p=0}^{n-1}a_pz^{-p}|<M_{Y,n}.$$ 

If there exist the constants $A_Y,\,C$ i.~e. $M_{Y,n}=C (n!)^{1/k}A_Y^n$ a series $\hat{w}$ is \textit{an asymptotical Gevrey series  of order $1/k$ for a function $w$  on $\Omega_R(\varphi_1,\varphi_2)$}. We define this like $\hat{w} \in \mathbb{C}[[z]]_{1/k}$.

For $G(z,Y,Y_1,\ldots,Y_{n})$ being an analityc function of $n+2$ variables. Let us consider an equation
\begin{equation}
\label{Gev}
G(z,w,Dw,\ldots,D^nw)=0.
\end{equation}
Let  $\hat{w}\in \mathbb{C}[[1/z]]$ be a formal series, being a formal solution of a differential equation (\ref{Gev}), and $D$ --- an operator $z\dfrac{d}{dz}$.




\begin{theorem}
\label{Gev-1} (O. Perron, J.-P. Ramis, B. Malgrange, Y. Sibuya in different cases)
Let $\hat{w}\in \mathbb{C}[[z]]_{1/k}$ be a solution to the equation $(\ref{Gev})$. Then there exist $k'>0$ i. e. for every open sector $V$ with the vertex in the infinity, having an angle less than $\min(\pi/k,\pi/k')$ and for a sufficiently big $R$ there exist a function  $w$, being a solution to the equation $(\ref{Gev})$ which is asymptotically approximated of Gevrey order $1/k$ by a series $\hat{w}$.
\end{theorem}
The next theorem contains conditions on the Newton polygon. We will describe the process of its construction.

Let us be given a linear differential operator
\begin{equation*}
L=\sum_{k=0}^{n}a_k(z)D^k,\,\mbox{where } a_k(z) \in \mathbb{C}[[z]][1/z],\,a_k(z)=\sum_{j_k=j_{k,\,0}}^\infty a_{j_k}z^{-j},$$$$ a_{j_k}=\mathrm{const}\in\mathbb{C}.
\end{equation*}
 We will put in correspondence a set of points on the plane: $(k,\,j_{k,\,0}),\,k=0,\,\ldots,\,n$ --- \textit{a support of the operator} L. We will define a set $$
\displaystyle{N=\bigcup_{k=0}^{n}\{(q_1,q_2):q_1\leq k,\,q_2\geq j_{k,\,0}\}}$$ and then we construct a convex hull of this set in half-plane $q_1\geq 0$. A boundary of this set is called \textit{a Newton polygon of the linear differential operator} $L$. 
\begin{theorem}
\label{Gev-2}(J.-P. Ramis) 
Let $\hat{w} \in \mathbb{C}[[z]]$ be a formal solution to the equation $(\ref{Gev})$, then the series $\hat{w}$ converges or has a Gevrey order equal to $s$, where $s\in \{0, \dfrac{1}{k_1},\,\ldots,\,\dfrac{1}{k_N}\}$ and  $0<k_1<\ldots<k_N<+\infty$ are all positive slopes of the edges of the $N(G,\hat{w})$ to the X-axis.
\end{theorem}

The previous theorem \ref{Gev-2} is formulated for a nonlinear differential equation and its formal solution  $\hat{w}$.

In a particular case $\dfrac{\partial G}{\partial Y_n}(z,\hat{w},\,,\ldots,\,\hat{w}^{(n)})\neq 0$ a Newton polygon of this equation on a formal solution can be constructed (Remark A.2.4.3  ~\cite{Sibuya}) as a polygon of an operator
\begin{equation}
\label{L}
L_0=\sum_{k=0}^n\left(\frac{\partial G}{\partial Y_k}(z,\hat{w},D\hat{w},\,\ldots,\,D^n \hat{w})\right)D^k.
\end{equation}
We can show that this operator coincides with the operator $\mathcal{M}$ used to construct exponential expansions of solutions using the methods of Power Geometry \cite{exp}.


Let $(a_{j,\,k})$ be a matrix of a transformation from the basis $D,\,D^2,\,\ldots,\,D^n$  to basis \newline$\dfrac{d}{dz},\, {z^2}\dfrac{d^2}{dz^2},\,\ldots,\, {z^n}\dfrac{d^n}{dz^n}$ in a vector space of linear differential operators. We can check that the element $a_{j,\,k}$   is equal to $a_{j,\,k}=S(j,k)$, where $S(j,k)$ is a Stirling number of the second kind.
This assertion is easily proved by induction using the fact that the elements of the matrix satisfy an equation $a_{j+1,\,k+1}=a_{j,\,k}+(k+1)a_{j,\,k+1}$. The transition matrix is lower triangular, the diagonal elements are equal to one ($S(j,j)=1$). An inverse matrix $(a^{j,k})$ is also lower triangular, its elements are equal to $a^{j,k}=(-1)^{j-k}s(j,k)$, where $s(j,k)$ is a Stirling number of the first kind.

Let $G(z,w,Dw,\,\ldots,\,D^n w)=F(z,w,z\dfrac{d}{dz},\,\ldots,\,z^n\dfrac{d^n}{dz^n}),$  i.~e.\newline $G(z,Y,Y_1,\ldots,\,Y_n)=F(z,Y,X_1,\,\ldots,\,X_n)$.
An operator (\ref{L}) can be rewritten using ${z^k}\dfrac{d^k}{dz^k}$ (we use designation $\delta_{j,l}$ for Kronecker delta): 
\begin{equation*}
L_0=\sum_{k=0}^n\sum_{l=0}^n\frac{\partial F}{\partial X_l}\frac{\partial X_l}{\partial Y_k}D^k=\sum_{k=0}^n\sum_{l=0}^n\frac{\partial F}{\partial X_l}a^{l,\,k}\sum_{j=0}^na_{k,\,j}z^j\dfrac{d^j}{dz^j}=$$ $$=\sum_{l=0}^n\dfrac{\partial F}{\partial X_l}\sum_{j=0}^nz^j\dfrac{d^j}{dz^j}\sum_{k=0}^na^{l,k}a_{k,j}=\sum_{l=0}^n\dfrac{\partial F}{\partial X_l}\sum_{j=0}^nz^j\delta_{j,l}\dfrac{d^j}{dz^j}=\sum_{l=0}^n\dfrac{\partial F}{\partial X_l}z^l\dfrac{d^l}{dz^l}, 
\end{equation*}
 and this expression is a first variation of a differential sum $F$. Being calculated on a series $\hat{w}$, it coincides with an operator $\mathcal{M}$ from \cite{exp}.

The conclusion is -- to construct a Newton polygon of the equation \newline$F(z,w,w',\ldots,w^{(n)})=0$ on a solution $\hat{w}$ we need to perform the following:

1) to calculate the first variation $\dfrac{\delta F}{\delta w}$ on a solution $\hat{w}$ ; 

2) to perform a transformation expressed by a matrix $(-1)^{j-k}s(j,k)$;

3) to verify a condition $\dfrac{\partial G}{\partial Y_n}(z,\hat{w},\,,D\hat{w},\,,\ldots,\,D^n\hat{w})\neq 0$; 

4) to find the set $(k,\,j_{k,\,0}),\,k=0,\,\ldots,\,n$;

5) to construct  a convex hull of this set in half-plane $q_1\geq 0$.

\textbf{Remark 1.} The steps 2 and 3 can be interchanged: i.~e. instead of the condition \newline$\dfrac{\partial G}{\partial Y_n}(z,\hat{w},\,,D\hat{w},\,,\ldots,\,D^n\hat{w})\neq 0$ we  can check the condition  $$\dfrac{\partial F}{\partial X_n}(z,\hat{w},\,,\hat{w}',\,,\ldots,\,\hat{w}^{(n)})\neq 0.$$


\textbf{Remark 2.} We use a Newton polygon considered in \cite{Sibuya} to calculate the Gevrey orders but it can be easily shown that Gevrey order of the solution can be calculated via constructing a polygon of the eqaution used in Power Geometry \cite{BrUmn}. To calculate the order we should perform a transrormation $u=\ln y$ in the operator $\mathcal{M}(z)$ \cite{exp}
, applied to $u$,  to reduce an expression by $u$ and for the differential sum obtained (the sum depends on $z$ and $y$) to construct a polygon. In conditions of the theorem \ref{Gev-2} the tangents are replaced by cotangents.

\section{The fifth Painlev\'{e} equation}
We consider the fifth Painlev\'{e} equation
\begin{equation}
\label{P5} w''=\left( \frac{1}{2w}+\frac{1}{w-1}
\right)\left(w'\right)^2-\frac{w'}{z}+\frac{(w-1)^2}{z^2}\left(\alpha
w+ \frac{\beta}{w}\right)+ \frac{\gamma w}{z}+\frac{\delta w
(w+1)}{w-1},
\end{equation}
where $\alpha, \beta, \gamma, \delta$
are complex parameters,  $z$ is an independent complex variable, $w$ is a dependent one, and we consider its formal power series solutions near infinity. Such solution are obtained in a work \cite{Gromak}. 

If $\alpha\beta\gamma \delta \neq 0$ there exist the following five expansions:
\begin{equation}
\label{solF_11} 
(-1)^l\sqrt{\frac{\beta}{\delta}}\frac{1}{z}
+\left(-\frac{2\beta}{\delta}+(-1)^l\frac{\gamma}{2\delta}\sqrt\frac{\beta}{\delta}\right)\frac{1}{z^2}+
\sum_{s=3}^{\infty} \frac{c_{-s,l}} {z^s},~
l=1, 2,
 \end{equation} 
\begin{equation}
\label{solF_2}  -1+\frac{2 \gamma}{\delta z}+
\sum_{s=2}^{\infty} \frac{c_{-s}} {z^s},
\end{equation}
\begin{equation}
\label{solF_31}
  (-1)^l
\sqrt{-~\frac{\delta}{\alpha}} z+ 2+(-1)^l\frac{\gamma}{2\sqrt{-\alpha \delta}}+
\sum_{s=1}^{\infty} \frac{c_{-s,l}} {z^s},~l= 1, 2.
\end{equation}
If $\alpha\beta\gamma \neq 0$, $\delta=0$ there exist the following four expansions:
\begin{equation}
\label{soldF_11}
(-1)^l\sqrt{-~\frac{\displaystyle \beta}{\displaystyle \gamma}}
\frac{\displaystyle 1}{\displaystyle \sqrt{z}}+ \frac{
\beta}{\gamma z}+\sum_{s=3}^{\infty} \frac{c_{-s,l}}
{z^{s/2}},~ l = 3, 4,
\end{equation}
\begin{equation}
\label{soldF_31}
(-1)^l
\sqrt{-\dfrac{\gamma}{\alpha}} \sqrt{z}+
1+ \sum_{s=1}^{\infty} \frac{c_{-s,l}} {z^{s/2}},~~l=3,4.
\end{equation}

The coefficients $c_s$, $c_{s,l} \in \mathbb{C}$, $l=1,2,3,4$  are uniquely determined constants, i.~e.: if we fix the values of the parameters $\alpha, \beta, \gamma, \delta$ these coefficients are uniquely determined as solutions of a non-degenerate system of linear equations.

\begin{theorem}
\label{Gevrey}
The series  $(\ref{solF_2})$, $(\ref{solF_31})$ and a regular part of a series $(\ref{solF_11})$ are of Gevrey order $1$. The series $(\ref{soldF_11})$, $(\ref{soldF_31})$ considered as the series in a new variable $\sqrt{z}$ are also of Gevrey order $1$.
\end{theorem}

\textbf{Proof.} We apply the theorem \ref{Gev-2} taking as an equation (\ref{Gev}) an equation (\ref{P5}) multiplied by $z^2 w (w-1)$ with all the terms of the equation put into the right part:
\begin{equation}
\label{P_mod}
 f(z,w) \stackrel{def}{=} -z^2 w (w-1) w''+z^2 \left(
 \frac{3}{2}w-\frac{1}{2}\right)\left(w'\right)^2
-z w (w-1) w'+$$
$$ +(w-1)^3(\alpha w^2+\beta)+
\gamma z w^2 (w-1)+\delta z^2 w^2 (w+1)=0,
\end{equation}
 we take instead  $\hat{w}$ the series  (\ref{solF_11}), (\ref{solF_2}), (\ref{solF_31}) in course. If the principal part is not equal to zero, we can easily obtain the case of a zero principal part using a transformation; we speal about Gevrey order of a regular prat of the series. 

The first variation of the equation $P_5$ represented in a form of a differential sum (\ref{P_mod}) is equal to 
\begin{equation}
\label{P_5_var}
-z^2w(w-1)\frac{d^2}{dz^2}+\left(z^2(3w-1)w'-zw(w-1)\right)\frac{d}{dz}-$$$$
-z^2(2w-1)w''+\frac{3z^2(w')^2}{2}-z(2w-1)w'+$$ $$+(w-1)^2(5\alpha w^2-2\alpha w+3\beta)+\gamma z(3w^2-2w)+\delta z^2(3w^2+2w).
\end{equation}

We substitute a series (\ref{solF_31}) to the expression (\ref{P_5_var}) and write only coefficients of $\dfrac{d^2}{dz^2}$, $\dfrac{d}{dz}$ and identity operator with the maximum degree in $z$:
\begin{equation*}
(-1)^l\sqrt{\dfrac{\beta}{\delta}}\left(z\frac{d^2}{dz^2}+2\frac{d}{dz}+2\delta z\right),\,l=1,2,
\end{equation*}
a support of such an operator consists of the points� $(0,-1),\,(1,1),\,(2,1),$ the Newton polygon is shown in Fig. 1.

Analogous calculations can be performed for the series (\ref{solF_2}). We obtain an operator
\begin{equation*}
-2z^2\frac{d^2}{dz^2}-z\frac{d}{dz}+\delta z^2,
\end{equation*}
the support of it consists of the points $(0,-2),\,(1,0),\,(2,0),$ its Newton polygon (brought down by a vector $(0,1)$) is shown in Fig. 1. 

For the series (\ref{solF_11}) we obtain an operator
\begin{equation*}
\dfrac{\delta}{\alpha}z^4\frac{d^2}{dz^2}+2 \sqrt{\dfrac{-\delta}{\alpha}}z^3\frac{d}{dz}-3\dfrac{\delta^2}{\alpha} z^4,\,l=1,\,2.
\end{equation*}
Its support consists of the points $(0,-4),\,(1,-2),\,(2,-2),$ its Newton polygon (brought down by a vector $(0,3)$) is shown in Fig. 1. 

As we see in Fig. 1, the unique positive tangent of the Newton polygon is equal to $1$, using the theorem \ref{Gev-2} we obtain that the series (\ref{solF_31}), the series (\ref{solF_2}) and a regular part of the series (\ref{solF_11}) are of Gevrey order $1$.  


Let us calculate Gevrey order of the series (\ref{soldF_11}) and (\ref{soldF_31}) obtained if $\delta =0$.
To transform the fifth Painlev\'{e} equation to the form (\ref{Gev}) we perform substitute $t=\sqrt{z}$ and consider the series (\ref{soldF_11}) and (\ref{soldF_31})  as the series in decreasing half-integer degrees of $t$. We calculate the first variation (we define differentiationg with the respect to  $t$ as a dot):
\begin{equation}
\label{P_5_var_delta}
-~\dfrac{t^2w(w-1)}{4}\frac{d^2}{dt^2}+\frac{1}{4}\left(t^2(3w-1)\dot{w}+tw(w-1)-w(w-1)\right)\frac{d}{dt}-$$$$
-~\frac{(2w-1)(t^2\ddot{w}+t\dot{w})}{4}+\frac{3t^2\dot{w}^2}{8}+$$ $$+(w-1)^2(5\alpha w^2-2\alpha w+3\beta)+\gamma t^2(3w^2-2w).
\end{equation}

We substitute a series  (\ref{soldF_31}) to the expression (\ref{P_5_var_delta}) and write only coefficients of $\dfrac{d^2}{dt^2}$, $\dfrac{d}{dt}$ and identity operator with the maximum degree in $t$:
\begin{equation*}
\dfrac{c_{-1/2}t}{4}\dfrac{d^2}{dt^2}+
\dfrac{c_{-1/2}-3c_{-1/2}^2}{4t}\dfrac{d}{dt}+3\gamma c_{-1/2}t=\dfrac{c_{-1/2}}{4t}D^2_t+\dfrac{-c_{-1/2}}{t}D_t+3\gamma c_{-1/2}t,
\end{equation*}
where $D_t=t\dfrac{d}{dt}$, $c_{-1/2}=(-1)^l\sqrt{\dfrac{\beta}{\gamma}},\, l=1,2$. 

 The support of the operator consists of the points $(0,-1),\,(1,1),\,(2,1),$ its Newton polygon is shown in Fig. 1. We obtain again that the series  (\ref{soldF_31}) is of Gevrey order $1$.

An operator corresponding to a series (\ref{soldF_11}) has the form
\begin{equation*}
c_{1/2}^2\dfrac{t^4}{4}\dfrac{d^2}{dt^2}+
c_{1/2}^2t^3\dfrac{d}{dt}-3\gamma c_{1/2}^2t^4=\dfrac{c_{1/2}^2t^2}{4}D^2_t-\dfrac{3c_{1/2}^2t^2}{4}D_t-3\gamma c_{1/2}^2t^4,
\end{equation*} 
where  $c_{1/2}=(-1)^l\sqrt{-\dfrac{\gamma}{\alpha}},\, l=1,2$.

The support of the operator consists of the points $(0,-4),\,(1,-2),\,(2,-2),$ its Newton polygon (brought down by a vector $(0,3)$) is shown in Fig. 1. The regular part of the series (\ref{soldF_11}) considered as the Laurent series in a variable $t$ is of Gevrey orser $1$.

This statement accomplishes the proof of the theorem \ref{Gevrey}. 

\textbf{Assertion 1.} \textit{
There exist  $k'\geq 1$ and $R_0\in \mathbb{R}_+$ i.~e. for every open sector  $\{z:|z|>R\geq R_0,\,\mathrm{Arg}\,z\in(\varphi_1,\varphi_2)\},$ $\varphi_2-\varphi_1<\pi/k'\leq\pi$ there exist a solution to the fifth Painlev\'{e} equation approximated by this Gevrey-1
series 
 by each of  the power series (series $(\ref{solF_31})$, series $(\ref{solF_2})$ and regular part of the series $(\ref{solF_11})$) obtained.}

\textit{There exist $k'\geq 1/2$ and $R_0\in \mathbb{R}_+$ i.~e. for every domain $\{z:|\sqrt{z}|>R\geq R_0,\,\mathrm{Arg}\,z\in(\varphi_1,\varphi_2)\},$  $\varphi_2-\varphi_1<\pi/k'\leq 2\pi$ there exist a solution to the fifth Painlev\'{e} equation (with a parameter $\delta=0$) approximated by these Gevrey-1
series  $(\ref{soldF_31})$ and  $(\ref{soldF_11})$ correspondingly}.

This assertion is consequence of the theorems \ref{Gev-1} and \ref{Gevrey}.



The results of this section are partially published as a preliminary version in \cite{Parus-per}.

Other results concerning the fifth Painlev\'{e} equation are published in \cite{Parus-neos}, \cite{Parus-loc}.

\section{The third Painlev\'{e} equation}

Let us pass on to the consideration of the analogous questions for the power expansions of oslutions to the third Painlev\'{e} equation:
\begin{equation}
 \label{P3}
w''=\frac{(w')^2}{w}-\frac{w'}{z}+\frac{\alpha
w^2+ \beta}{z}+ \gamma w^3+\frac{\delta}{w},
\end{equation}
which can also be found in a book \cite{Gromak}.

When $\alpha\beta\gamma\delta\neq 0$ they form four power expansions:
\begin{equation}
 \label{solP3}i^l\sqrt[4\,]{-\dfrac{\delta}{\gamma}}-\left(\dfrac{(-1)^l\beta}{4\sqrt{-\gamma\delta}}+\dfrac{\alpha}{4\gamma}\right)\dfrac{1}{z}+\sum_{k=2}^\infty \dfrac{c_{l,k}}{z^k},\,l=1,2,3,4,
\end{equation}
where $i^2=-1,$ and we consider the main branch of the root while speaking about the forth root.

We apply the theorem \ref{Gev-2} taking  an equation (\ref{P3}) multiplied by $z w $ with all the terms of the equation put into the right part as an equation (\ref{Gev}):
\begin{equation}
\label{P3_mod}
  -z w  w''+z \left(w'\right)^2
- w  w'+w(\alpha w^2+\beta)+
\gamma z w^4 +\delta z=0,
\end{equation}
 and taking the series (\ref{solP3}) as a formal power series solution $\hat{w}$. 

The first variation of the equation  (\ref{P3_mod}) is equal to 
\begin{equation}
\label{P_3_var}
-zw\frac{d^2}{dz^2}+\left(2zw'-w\right)\frac{d}{dz}-zw''-w'+3\alpha w^2+\beta+4\gamma z w^3.\end{equation}

We substitute a series (\ref{solP3}) to the expression (\ref{P_3_var}) and write only coefficients of $\dfrac{d^2}{dz^2}$, $\dfrac{d}{dz}$ and identity operator with the maximum degree in $z$:
\begin{equation*}
i^k\sqrt[4\,]{-\dfrac{\delta}{\gamma}}z\frac{d^2}{dz^2}+B\frac{d}{dz}+4 zi^{3k}\sqrt[4\,]{{\delta^3}{\gamma}},\,l=1,2,3,4,
\end{equation*}
a support of such an operator consists of the points� $(0,-1),\,(1,1),\,(2,1),$ the Newton polygon is shown in Fig. 1.

As we see in Fig. 1, the unique positive tangent of the Newton polygon is equal to $1$, using the theorem \ref{Gev-2} we obtain that the series  (\ref{solP3})  are of Gevrey order $1$. 

So, we obtain the followin theorem
\begin{theorem}
\label{GevreyP3}
The series  $(\ref{solP3})$ are of Gevrey order equal to one. There exist $k'\geq 1$ � $R_0\in \mathbb{R}_+$ i.~e. for any open sector $\{z:|z|>R\geq R_0,\,\mathrm{Arg}\,z\in(\varphi_1,\varphi_2)\},$ $\varphi_2-\varphi_1<\pi/k'\leq\pi$ there exist a solution to the third Painlev\'{e} equation approximated by this Gevrey-1
series.
\end{theorem}

\begin{picture}(70, 60)(0,50)
\put(-10,0){\vector(1,0){80}}
 \put(0,-50){\vector(0,1){140}}
\put(-4,20){\line(1,0){8}} \put(-4,40){\line(1,0){8}}\put(-4,-20){\line(1,0){8}}
 {\Large \put(60,-10){$q_1$}
 \put(-10,80){$q_2$}
 \put(-10,-10){$0$}
\put(-20,-25){$-1$}
 \put(38,-12){$2$}
 \put(15,-75){Fig. 1}
  \thicklines{
    \put(0,-40){\line(1,1){40}}
    \put(0,-40){\line(0,1){100}}
    \put(40,0){\line(0,1){80}}
    }
\put(0,-40){\circle*{3}}
  \put(20,0){\circle*{3}}
\put(40,0){\circle*{3}}}
\end{picture}

\vskip 2 cm
National Research University Higher School of Economics,\\
Moscow Institute of Electronics and Mathematics,\\
Bolshoi Trekhsvjatitelskii per. 3, Moscow, 109028, Russia\\
e-mail: parus-a@mail.ru

\end{document}